\documentclass[11pt]{amsart}
 \usepackage{amsaddr}
 
\usepackage[T1]{fontenc}
\usepackage{mathpazo}
\usepackage{geometry}
\usepackage{amsmath}
\usepackage{amsfonts}
\usepackage{amssymb}
\usepackage{amsthm}
\usepackage{mathtools}
\usepackage{hyperref}
\usepackage{tikz}
\usepackage{graphicx}
\usepackage{pgfplots}
\usepackage{subcaption}

\makeatletter
\@namedef{subjclassname@2020}{%
	\textup{2020} Mathematics Subject Classification}
\makeatother

\newtheorem{theorem}{Theorem}[section]
\newtheorem*{theorem*}{Theorem}
\newtheorem{proposition}[theorem]{Proposition}
\newtheorem*{proposition*}{Proposition}
\newtheorem{corollary}[theorem]{Corollary}
\newtheorem*{corollary*}{Corollary}
\newtheorem{lemma}[theorem]{Lemma}
\newtheorem*{lemma*}{Lemma}
\newtheorem{conj}{Conjecture}
\newtheorem*{conj*}{Conjecture}

\theoremstyle{definition}

\newtheorem{remark}[theorem]{Remark}
\newtheorem{example}[theorem]{Example}

\renewcommand{\tilde}{\widetilde}
\renewcommand\epsilon{\varepsilon}

\newcommand\F{\mathbb{F}}
\newcommand\Fq{\mathbb{F}_q}

\newcommand\PP{\mathbb{P}}
\newcommand{\CPhi}{C_\Phi}

\title{
	Curves on Frobenius classical surfaces in $\PP^3$ over finite fields
}

\author{Elena Berardini}
\address{LTCI, Télécom Paris, Institut polytechnique de Paris, F-91120 Palaiseau, France}
\email{elena\_berardini@hotmail.it}

\author{Jade Nardi}
\address{Univ Rennes, CNRS, IRMAR - UMR 6625, F-35000 Rennes, France}
\email{jade.nardi@univ-rennes1.fr}


\begin{document}
\sloppy
	\maketitle
	
	\begin{abstract} 
	In this paper we give an upper bound on the number of rational points on an irreducible curve $C$ of degree $\delta$ defined over a finite field $\F_q$ lying on a Frobenius classical surface $S$ embedded in $\PP^3$. This leads us to investigate arithmetic properties of curves lying on surfaces. 
	In a certain range of $\delta$ and $q$, our result improves all other known bounds in the context of space curves. 
	\end{abstract}

\section*{Introduction}
The number of rational points on a smooth projective absolutely irreducible curve $C$ of genus $g$ defined over the finite field $\F_q$ is bounded by the famous Serre--Weil bound, namely $\#C(\F_q)\leq q+1+g\lfloor 2\sqrt{q}\rfloor$. This bound was extended to possibly singular curves by Aubry and Perret \cite{AP11}. In their pivotal work of 1986, Stöhr and Voloch \cite{SV86} introduced the Frobenius orders of a projective smooth curve and used them to give an upper bound on the number of rational points on the curve. Nearly 30 years later, Homma proved that the number of $\F_q$--points on a non--degenerate curve of degree $\delta$ embedded in $\PP^n$, with $n\geq 3$, does not exceed $q(\delta-1)+1$  \cite{H12}. 
All these results improve the original Serre--Weil bound for a range of parameters, and often deal with more general curves (possibly reducible and/or singular). Such bounds are interesting in themselves, and also turn out to be useful for applications to coding theory. Indeed, when studying algebraic geometry codes from surfaces, one approach to bound the minimum distance of such codes relies on bounding the number of rational points on irreducible curves that appear as components of some divisors on the surface (see for instance \cite{VZ09,MYFE}). The sharper the upper bounds for the number of $\F_q$--points on those curves, the sharper the lower bounds for the minimum distance.

\vspace{5pt}
In this paper, we deal with irreducible curves lying on surfaces embedded in $\PP^3$. Our method is inspired by the techniques developed by Stöhr and Voloch in \cite{SV86}. More precisely, let $S$ be an irreducible algebraic surface defined over a finite field $\Fq$ and embedded in a projective space $\PP^3$. Consider the set $\CPhi^S$ of points $P$ on $S$ whose image under the Frobenius map lies on their own tangent plane, \textit{i.e.} $\CPhi^S  \coloneqq \{P \in S \mid \Phi(P) \in T_PS\}$. When the variety  $\CPhi^S$ is a curve, the surface is called Frobenius classical.

Assume that $S$ is a Frobenius classical surface and take $C$ a curve lying on $S$. Obviously, the $\Fq$--points on $C$ are contained in the intersection of $C$ with $\CPhi^S$.  If $C \cap \CPhi^S$ consists of a finite set of points, we can bound the number of rational points on $C$ by the cardinality of $C \cap \CPhi^S$. A classical way to sharpen this bound is to prove that the $\F_q$--points appear with multiplicity greater than $1$ in this intersection (see \cite{SV86, V03}). This approach leads to the following result. 

\begin{theorem*}\label{th:intro}(Theorem~\ref{th:P3})
	Let $S$ be an irreducible Frobenius classical surface defined over $\Fq$ of degree $d>1$ in $\PP^3$. Let $C$ be an irreducible curve defined over $\Fq$ of degree $\delta$ lying on $S$. Suppose $C$ is not a component of $\CPhi^S$. Then
	\begin{equation}\tag{\ref{notreborne}}
		\#C(\F_q)\leq \frac{\delta(d+q-1)}{2}. 
	\end{equation}

\end{theorem*}

 If the curve $C$ shares some components with $\CPhi^S$, our strategy fails. A large part of the present work is then dedicated to understanding the components of the curve $\CPhi^S$. With this purpose in mind, we distinguish Frobenius classical and non--classical curves, and determine their positioning with respect to $\CPhi^S$ depending on their Frobenius orders $\nu_1$ and $\nu_2$ (as introduced in \cite{SV86} and recalled in Section \ref{Sec:geom_of_curves}). We obtain the following results.
 
 \begin{proposition*}[Propositions \ref{prop:FNCA=Nice} and Proposition \ref{cor:nu1>1}]
Non--plane Frobenius non--classical curves are components of $\CPhi^S$ if and only if $\nu_1>1$. Moreover, curves with $\nu_1>1$ of degree at most $q$ are plane curves.
\end{proposition*} 	

Regarding Frobenius classical components, we conjecture the following feature of $\CPhi^S$.

\begin{conj*}[Conjecture \ref{conjecture}]
Non--plane Frobenius classical components of $\CPhi^S$ have degree larger than $q$.
\end{conj*}
%
 
 This conjecture is guided by several experiments performed with the {\scshape{Magma}} computer algebra system, but we are not able to prove it. However, it seems an interesting question in itself, and it only affects our result in the case of Frobenius classical curves, when they are components of $\CPhi^S$.
 
 
%
%
%

\medskip

The paper is organised as follows. In Section \ref{Sec:geom_of_curves}, we recall some general results on the theory of orders of a space curve. Theorem \ref{th:P3} is proved in Section \ref{Sec:borne}. Section \ref{Sec:components} gathers results on curves lying on a surface embedded in $\PP^3$ to handle the components of $\CPhi^S$. We also state and discuss Conjecture \ref{conjecture}. Finally, in Section \ref{Sec:conclusion}, we compare our bound to existing ones and present future research directions.

\section{Geometry of space curves}\label{Sec:geom_of_curves}

Let us recall the notion of order at a point and Frobenius orders for space curves, and some of their properties. We refer the reader to \cite{HKT13,SV86} for a more general discussion and the proofs of the results stated here.

\medskip
Let $\Fq$ be a finite field of size $q$ and let $\overline{\F}_q$ denote an algebraic closure of it. 
Let $C \subset \PP^3$ be an absolutely irreducible projective  space curve defined over $\overline{\F}_q$.  Let $P\in C$. An integer $j$ is called a \emph{$P$--order} if there exists a plane intersecting the curve $C$ with multiplicity $j$ at $P$. By definition of the degree of a curve, the $P$--orders $j_i$ are bounded from above by $\deg C$.

Assume that $C$ is \emph{non--degenerate}, that is, not contained in a plane. If $P$ is a non--singular point of $C$, there are exactly four distinct $P$--orders $j_0=0 < j_1< j_2 < j_3$ \cite[Theorem~1.1]{SV86}. Moreover, in this case, the first nontrivial $P$--order $j_1$ equals $1$.

We denote by $T_P^{(i)}C$ the intersection of the hyperplanes meeting the curve $C$ at $P$ with multiplicity at least $j_{i+1}$. Then 
\[P=T_P^{(0)}C\subset T_P^{(1)}C \subset T_P^{(2)}C \subset \PP^3.\]
If $P$ is non--singular, the space $T_P^{(1)}C$ is the tangent line. The projective plane $T_P^{(2)}C$ is called the \emph{osculating plane} of $C$ at $P$.

Let $t$ be a local parameter at $P$. The Hasse derivatives $D_t^{(i)}$ with respect to $t$ form a family of $\overline{\F}_q$--linear operators ($i \geq 0$) on the function field $\overline{\F}_q(C)$ such that 
\[D_t^{(i)}t^k=\binom{k}{i}t^{k-i} \text{ and } D_t^{(i)}D_t^{(k)}=\binom{i+k}{i}D_t^{(i+k)} \text{ for } i,k \geq 0.\]
The $P$--orders give the equation of the osculating plane \cite[Corollary~1.3]{SV86}~: 
\[T_P^{(2)}C:\det\begin{pmatrix}
	X_0 			& X_1			 	& X_2 					& X_3 \\
	x_0 			&x_1 				& x_2					& x_3 \\
	D_t^{(j_1)}x_0 	& D_t^{(j_1)}x_1 	& D_t^{(j_1)}x_2		& D_t^{(j_1)}x_3 \\
	D_t^{(j_2)}x_0 	& D_t^{(j_2)}x_1 	& D_t^{(j_2)}x_2		& D_t^{(j_2)}x_3 \\
\end{pmatrix}=0.\]
For almost every point $P \in C$, the sequence of $P$--orders is the same, say $(\epsilon_0,\epsilon_1,\epsilon_2,\epsilon_3)$. There are only finitely many points such that
$(j_0,j_1,j_2,j_3)\neq(\epsilon_0,\epsilon_1,\epsilon_2,\epsilon_3)$, which are called the \emph{Weierstrass points} of the curve. For a curve of genus $g$, the number of such points, counted with multiplicity, equals $(\epsilon_1+\epsilon_2+\epsilon_3)(2g-2)+4\deg(C)$ \cite[\S 1]{SV86}.

One always has $\epsilon_0=0$ and $\epsilon_1=1$, as a curve is non--singular at almost all of its points. A curve is said to be \emph{classical} if $(\epsilon_0,\epsilon_1,\epsilon_2,\epsilon_3)=(0,1,2,3)$ and \emph{non--classical} otherwise. Homma gave the list of all the possible orders sequence of space curves \cite{H91}.
 The geometry of non--classical curves was studied by Hefez and Kakuta \cite{HK92}. 

Take $P \in C$. Its image under the Frobenius map belongs to its osculating plane if and only if 
\[\det\begin{pmatrix}
	x_0^q 			& x_1^q			 	& x_2^q 				& x_3^q \\
	x_0 			&x_1 				& x_2					& x_3 \\
	D_t^{(j_1)}x_0 	& D_t^{(j_1)}x_1 	& D_t^{(j_1)}x_2		& D_t^{(j_1)}x_3 \\
	D_t^{(j_2)}x_0 	& D_t^{(j_2)}x_1 	& D_t^{(j_2)}x_2		& D_t^{(j_2)}x_3 \\
\end{pmatrix}=0\]
This motivates the introduction of the Frobenius order sequence of a curve, as follows.
\begin{proposition}\cite[Proposition~2.1]{SV86}\label{propSV}
	There exist integers $\nu_1 < \nu_2$ such that 
\[\det\begin{pmatrix}
	x_0^q 			& x_1^q			 	& x_2^q 				& x_3^q \\
	x_0 			&x_1 				& x_2					& x_3 \\
	D_t^{(\nu_1)}x_0 	& D_t^{(\nu_1)}x_1 	& D_t^{(\nu_1)}x_2	& D_t^{(\nu_1)}x_3 \\
	D_t^{(\nu_2)}x_0 	& D_t^{(\nu_2)}x_1 	& D_t^{(\nu_2)}x_2	& D_t^{(\nu_2)}x_3 \\
\end{pmatrix}\]
is a nonzero function. Choose them minimally with respect to the lexicographic order. Then
	\[\{\nu_1,\nu_2\} \subset \{\epsilon_1,\epsilon_2,\epsilon_3\}.\]
	
\end{proposition}

The integers $\nu_0=0,\nu_1,\nu_2$ are called the \emph{Frobenius orders} of $C$. The deleted order among $\{\epsilon_1,\epsilon_2,\epsilon_3\}$, say $\epsilon_I$, is called the \emph{$q$--deleted order} and the index $I$ is called the \emph{$q$--Frobenius index} of the curve. Garcia and Homma \cite{GH} studied both arithmetical and geometrical conditions satisfied by the $q$--deleted order of a curve. 
%
%
An interesting property of the $q$--deleted order can be deduced from Proposition \ref{propSV} :  for $i \in \{1,2\}$, 
\begin{equation}\label{eq:deleted_order}
I\leq i \: \Leftrightarrow \: \forall \: P \in C, \: \Phi(P)\in T_P^{(i)} C.
\end{equation}


Curves with Frobenius orders $\nu_i=i$ for every $i$ are called \emph{Frobenius classical}, otherwise they are called \emph{Frobenius non--classical}. 
Smooth Frobenius non--classical curves with $\nu_1 > 1$ have been extensively studied by Hefez and Voloch \cite{HV90}, who established their exact number of rational points. In the case of plane curves, singular Frobenius non--classical curves were investigated by Borges and Homma \cite{BH17}. 

A Frobenius non--classical curve is in general non--classical too. Examples of Frobenius non--classical curves that are classical are rare but do exist \cite[Remark 8.52]{HKT13}. However, for space curves, such examples can only occur when the $q$--Frobenius index $I$ satisfies $I\in\{1,2\}$ and $I+1 \equiv 0 \mod p$. Thus, in characteristic different from $2$ and $3$, Frobenius non--classical curves are necessarily non--classical too.

	\section{Rational points on curves over Frobenius classical surfaces}\label{Sec:borne}
	This section is devoted to the proof of Theorem \ref{th:P3}. To start with, we introduce Frobenius classical surfaces.
	\medskip
	
		Let $S \subset \PP^3$ be an irreducible projective normal surface defined by $f=0$ of degree $d$.  Set
	\begin{equation}\label{eq:h}
		h:=X_0^q f_0+X_1^q f_1+X_2^q f_2 + X_3^qf_3,
	\end{equation}
	where $f_i$ denotes the partial derivative of $f$ with respect to the $i$--th coordinate. 
	We define the variety $\CPhi^S$ as the zero set $f=h=0$. Set--theoretically, $\CPhi^S$ is the union of the singular locus of $S$ and the set of non--singular points $P \in S$ such that $\Phi(P) \in T_PS$. Note that since $S$ is normal, its singularities are isolated points \cite[Proposition 8.23]{H77}.\\
	A surface $S$ is said to be \emph{Frobenius classical} if $\CPhi^S$ is a curve, \textit{i.e.}, if $\dim \CPhi^S = 1$. In this case, since the polynomials $f$ and $h$ are of degree $d$ and $d+q-1$ respectively, the curve $\CPhi^S$ has degree $d(d+q-1)$.  Otherwise, the variety $\CPhi^S$ is the whole surface $S$ and we say that $S$ is \emph{Frobenius non--classical} (see \cite[Definition 2.1]{Asgarli} for an equivalent definition). For instance, the Hermitian surface defined over $\F_{q}$ by $X_0^{\sqrt{q}+1}+X_1^{\sqrt{q}+1}+X_2^{\sqrt{q}+1}+X_3^{\sqrt{q}+1}=0$ is Frobenius non--classical. 
	
	We have the following condition for a surface of degree $d$ to be Frobenius classical.
	\begin{lemma}\label{FCS}
	Let $S$ be an irreducible surface of degree $d$ defined over a finite field of characteristic $p$. Then, if $p\nmid d(d-1)$, $S$ is Frobenius classical.
	\end{lemma}
	\begin{proof}
	The hypothesis $p\nmid d(d-1)$ guarantees that the determinant of the Hessian matrix of $S$ does not vanish identically on $S$ \cite[Proposition 1]{V03}. This  implies that the surface $S$ is Frobenius classical by \cite[Theorem 4.6]{Asgarli}.
	\end{proof}
	Let us highlight that Lemma \ref{FCS} gives only a sufficient condition, as illustrated by the following example. 
	
	\begin{example}\label{ex:FCS}
		The absolutely irreducible surface defined by $f=X_0^3 + X_1^3 + X_2^2X_3=0$ over $\F_4$ is Frobenius classical. Indeed, in this case $h=X_0^6 + X_1^6 + X_2^2X_3^4=f^2+X_2^2X_3^2(X_3^2-X_2^2)$ is not a multiple of $f$.
	\end{example} 

	To bound the number of $\Fq$--points on a curve $C \subset S$, we shall count the number of points in the intersection $C \cap \CPhi^S$.

	\subsection{Intersection multiplicity of a curve $C$ and $\CPhi^S$ at rational points}\label{subsec:int_mult}

 	In order to sharpen the aforementioned bound, we will compute the intersection multiplicity of $C$ and $\CPhi^S$ at the rational points of $C$. More precisely, we will prove that the tangent spaces of the surface $S$ and the surface defined by $h=0$ (see Equation \eqref{eq:h}) at a non--singular rational point of $S$ are the same (Lemma \ref{lem:same_TS}). This will imply that the intersection multiplicity we consider is at least $2$ (Lemma \ref{lem:sing}).
 		  		 
	\begin{lemma}\label{lem:same_TS}
		Let $P$ be a non--singular rational point of the surface $S$ defined by $f=0$ of degree $d$. Let $S'$ be the surface defined by $h=0$, where $h$ is given in Equation \eqref{eq:h}. Then either  $T_P S = T_P S'$ or $P$ is a singular point of $S'$.\end{lemma}
	
	\begin{proof}
		Let us compare the gradients of $f$ and $h$ at $P$.
		The derivative of $h$ with respect to the $i$--th coordinate is equal to
		\[h_i =  X_0^q f_{i0}+X_1^q f_{i1}+X_2^q f_{i2}+X_3^q f_{i3}.\]
		Evaluating at $P$ and using Euler Lemma (see \cite[Exercise 5.8]{H77} for instance), we get $h_i(P)=(d - 1) f_i(P)$, where $d$ is the degree of $S$. This means that if the characteristic $p$ of the field does not divide $d-1$, then $T_P S = T_P S'$. Otherwise, the point $P$ is singular on $S'$.
	\end{proof}
		
	\begin{lemma}\label{lem:sing}
		Let $C \subset S$ be a curve defined over $\Fq$. Then for every $P \in C(\Fq)$ the intersection multiplicity of the curves $C$ and $\CPhi^S$ at $P$ is at least $2$.
	\end{lemma}

\begin{proof}
	Let us use the same notations as in the previous lemma. The curve $\CPhi^S$ is the complete intersection defined by $f=h=0$.  Then the intersection multiplicity of $C$ and $\CPhi^S$ at a point $P$ agrees with the intersection multiplicity of the curve $C$ and the surface $S'$ in $\PP^3$, which we denote by $i(P,C\cdot S';\PP^3)$. Applying \cite[Corollary 12.4]{Fulton} in $\PP^3$, we know that 
 	\[ i(P,C\cdot S';\PP^3) \geq m_P(C) m_P(S'),\]
 	where $m_P(C)$ and $m_P(S')$ are the multiplicities of $P$ in the varieties $C$ and $S'$, respectively. Moreover, the equality holds if and only if the projectivized tangent cones of $C$ and $S'$ at $P$ do not meet. 	
 	
 	If $P$ is a singular point on $C$ or on $S$', then $i(P,C\cdot S';\PP^3) \geq \max(m_P(C),m_P(S')) \geq 2$.
 	
 	Otherwise, $P \in C(\Fq) \subset S(\Fq)$ and $T_P S= T_P S'$ by Lemma \ref{lem:same_TS}. It means that the tangent plane of $S'$ at $P$ contains the tangent line of the curve $C$ at $P$, whence a non trivial intersection of the projectivized tangent cones. Then $i(P,C\cdot S';\PP^3) > 1$.
	
	In both cases we have $i(P,C\cdot S';\PP^3) \geq 2$, which concludes the proof.
\end{proof}

\subsection{Main result}
	
	\begin{theorem}\label{th:P3}
		Let $S$ be an irreducible Frobenius classical surface defined over $\Fq$ of degree $d>1$ in $\PP^3$. Let $C$ be an irreducible curve defined over $\Fq$ of degree $ \delta$ lying on $S$. Suppose $C$ is not a component of $\CPhi^S$. Then
		\begin{equation}\label{notreborne}
		\#C(\F_q)\leq \frac{\delta(d+q-1)}{2}. 
		\end{equation}
	\end{theorem}

	\begin{proof}
	Since $S$ is Frobenius classical, the set $\CPhi^S =\{P \in S \mid \Phi(P) \in T_PS\}$ is a curve. Clearly, we have ${C(\F_q) \subset C \cap \CPhi^S=\{P \in C \mid \Phi(P) \in T_PS\}}$. Since $C\not\subset \CPhi^S$, the curves $C$ and $\CPhi^S$ have no common components. 	
	By Lemma \ref{lem:sing}, the intersection multiplicity of $C$ and $\CPhi^S$ at every rational point on $C$ is larger than or equal to $2$. The space version of Bézout's theorem \cite[Chapter I, Theorem 7.7]{H77} applied on the curve $C$ and the surface defined by $h=0$ (see Equation \eqref{eq:h}) implies that the set $\CPhi^S \cap C$ has degree $\delta(d+q-1)$, which proves the result.
	\end{proof}

The bound in Theorem \ref{th:P3} is not necessarily an integer, so a trivial improvement can be obtained by taking the integral part of it.
\begin{remark}
	Even though our result is only stated for non--degenerate curves, our bound also applies to plane curves lying on a surface embedded in $\PP^3$, if they have at least one point whose image under the Frobenius map does not lie on its tangent line. This follows easily from \cite[Theorem 0.1]{SV86} using the fact that a plane curve on a degree $d$ surface has degree $\delta \leq d$.
	\end{remark}

The curve $\CPhi^S$ being defined over $\Fq$, if an irreducible curve $C$ that is not absolutely irreducible shares some components with $\CPhi^S$, then it is contained in $\CPhi^S$. Therefore, in order to handle the intersection with the curve $\CPhi^S$, we can focus on absolutely irreducible curves only and thus use the notions introduced in Section \ref{Sec:geom_of_curves}.

	\section{Properties of curves on space surfaces}\label{Sec:components}
	This section aims at characterising the components of the curve $\CPhi^S$, when $S$ is a Frobenius classical surface.  
	Our main goal here is to describe for which curves $C$, $C$ is not a component of $\CPhi^S$ and thus satisfies the main hypothesis of Theorem \ref{th:P3}. However, as a side effect, we also obtain some properties of curves lying over Frobenius non--classical surfaces.
	Before entering into details, let us observe the behavior of the Frobenius map under some change of coordinates.
	
	\subsection{Change of coordinates and Frobenius map}\label{subsec:coord_frob}
	
	Let $C \subset \PP^3$ be a curve defined over $\overline{\F}_q$. When studying local geometric properties around a point $P$ on $C$, one traditionally changes coordinates to make $P$ the origin. Before changing coordinates, the $q$--Frobenius map sends a tuple $(x_0:x_1:x_2:x_3)$ to its $q$--th power $(x_0^q:x_1^q:x_2^q:x_3^q)$. However, if $P$ is not a rational point, the Frobenius map after changing coordinates is no longer the $q$--th power.

	Let us fix a non--singular point $P=(x_0(P):x_1(P):x_2(P):x_3(P))$ of $C$. Let $t$ be a local parameter at $P$. Denote by $j_i$ the $P$--orders of $P$, for $i\in\{0,1,2,3\}$. As in \cite[Theorem 7.47]{HKT13}, there exists an invertible $4 \times 4$ matrix $A=(a_{ij})$ with coefficients in $\overline{\F}_q$ such that
	\[a_{i0}x_0+a_{i1}x_1+a_{i2}x_2+a_{i3}x_3=t^{j_i}+\dots\]

	Note that the point $P$ is rational if and only if $A$ has coefficients in $\Fq$.

	Setting new coordinates $(y_0,y_1,y_2,y_3)$ as $y_i=a_{i0}x_0+a_{i1}x_1+a_{i2}x_2+a_{i3}x_3$
	maps the point $P$ to the origin $(1:0:0:0)$ and $y_i=t^{j_i}+\dots$.
	
	Let us denote by $(\tilde y_0,\tilde y_1,\tilde y_2,\tilde y_3)$ the image of $(y_0,y_1,y_2,y_3)$ under the Frobenius map, \textit{i.e.} $\tilde y_i=a_{i0}x_0^q+a_{i1}x_1^q+a_{i2}x_2^q+a_{i3}x_3^q$. In particular, the Hasse derivatives $D^{(j)}_t \tilde{y_i}$ are zero for $j < q$, even when the point $P$ is not $\Fq$--rational.
	
%
%
	%

	\subsection{Components of $\CPhi^S$}\label{sec:compCPhi}
	The goal of this subsection is to get a better grasp on the components of $\CPhi^S$ and to determine some conditions implying that  a curve $C$ is not among them. It is clear that if $S$ contains some lines defined over $\F_q$, these lines are components of $\CPhi^S$. To address the problem for general curves, we will treat separately Frobenius non--classical curves with $\nu_1=1$ and $\nu_2>2$ (Section \ref{ss:FNC1}) and with $\nu_1>1$ (Section \ref{ss:FNC}), and Frobenius classical curves (Section \ref{ss:FC}). 
	Before diving into the proofs, let us set some notation.
	
	Fix a curve	$C$ lying on a normal Frobenius classical surface $S$. When dealing with a generic point $P$ on $C$, we can assume without loss of generality that $P$ is a non--singular point on both the surface $S$ and the curve $C$. We choose affine coordinates such that $P=(0,0,0)$ and $S$ and $C$ are locally given by 
	\begin{align}&S: z=u(x,y),&  C :\begin{cases}
		y=g(x), \\
		z=u(x,g(x)),
	\end{cases}
	\end{align}
	for some power series $u$ and $g$. In this system of coordinates, we denote by $(\tilde{x},\tilde{y},\tilde{z})$ the image under the Frobenius map of a point $(x,y,z)$.
	
	To ease the reading, we denote by $u^{(j)}$ and $g^{(j)}$ the $j$--th derivative of $u(x,g(x))$ and $g(x)$ with respect to $x$ for every $j \geq 0$. 
	
	For integers $1 \leq i < j $, we define the function $\Delta(i,j)$ by
	\begin{equation}\label{eq:Delta_ij}
		 \Delta(1,j)\coloneqq\det\begin{pmatrix}
			1			& \tilde{x}			& \tilde{y} 	& \tilde{z} \\
			1 & x		& y & z \\
			0 & 1 	& g'	& u'\\
			0 & 0 & g^{(j)}	& u^{(j)}
		\end{pmatrix}
		\text{ and }
		 \Delta(i,j)\coloneqq\det\begin{pmatrix}
			1			& \tilde{x}			& \tilde{y} 	& \tilde{z} \\
			1 & x		& y & z \\
			0 & 0 	& g^{(i)}	& u^{(i)}\\
			0 & 0 & g^{(j)}	& u^{(j)}
		\end{pmatrix} \text{ for } i \geq 2.\end{equation}
We recall that Proposition \ref{propSV} ensures the existence of two integers $\nu_1$ and $\nu_2$ such that $\Delta(\nu_1,\nu_2)$ is a nonzero function for a non--degenerate curve. Furthermore $\nu_1$ and $\nu_2$ are taken to be the minimal integers for which $\Delta(\nu_1,\nu_2)$ is a nonzero function. 
This proposition enables us to prove the following useful lemma that we will largely use in the rest of the section.

\begin{lemma}\label{lem:delta}
	Assume that we have $u^{(j)}=g^{(j)}u_y \text{ for every }j \geq \max\{2,\nu_1\}$, where $\nu_1$ denotes the first Frobenius order of $C$. Then either $\nu_1>1$ and $C$ is degenerate, or $\nu_1=1$ and, if $C$ is non--degenerate, then $\Phi(P) \notin T_P S$ for a generic point $P \in C$.
\end{lemma}
\begin{proof}
	We fix a generic point $P \in C$ and use the notations above.
	
		Assume $\nu_1 > 1$. Since for $j\geq \nu_1$ we have $u^{(j)}=g^{(j)}u_y$, we obtain
	\[\Delta(\nu_1,j)=\det\begin{pmatrix}
		1			& \tilde{x}			& \tilde{y} 	& \tilde{z} \\
		1 & x		& y & z \\
		0 & 0 	& g^{(\nu_1)}	& g^{(\nu_1)}u_y\\
		0 & 0 & g^{(j)}	& g^{(j)} u_y
	\end{pmatrix}=0,\]
	as the two last rows are collinear. This contradicts Proposition \ref{propSV} that ensures the existence of an integer $\nu_2$ such that $\Delta(\nu_1,\nu_2)\neq 0$. Consequently, the curve must be degenerate.\\	
	%
	Assume now that $\nu_1=1$. Using  that $u^{(j)}=g^{(j)}u_y$ for $j\geq 2$ we get
	\[\Delta(1,j)=g^{(j)}\det\begin{pmatrix}
		1			& \tilde{x}			& \tilde{y} 	& \tilde{z} \\
		1 & x		& y & z \\
		0 & 1 	& g' & u_x+g'u_y\\
		0 & 0 & 1	& u_y
	\end{pmatrix}=g^{(j)}\left[(\tilde{x}-x)u_x+(\tilde{y}-y)u_y-(\tilde{z}-z)\right].\]
	Since we have $\nu_1=1$ and $\Delta(1,\nu_2)\neq 0$ if $C$ is non--degenerate (by Proposition \ref{propSV}), the second factor must be nonzero, hence $\Phi(P) \notin T_P S$.
\end{proof}

\subsubsection{Frobenius non--classical curves with $\nu_1=1$.}\label{ss:FNC1}
We are going to prove that non--degenerate curves with $\nu_1=1$ and $\nu_2>2$ cannot be components of $\CPhi^S$.
\begin{proposition}\label{prop:FNCA=Nice}
		Assume that $C$ is non--degenerate and Frobenius non--classical with $\nu_1=1$. Then, for a generic point $P \in C$, we have $\Phi(P) \notin T_P S$. 
	\end{proposition}
		
	\begin{proof}
		Assume by contradiction that $\Phi(P) \in T_P S$. Since $C$ is Frobenius non--classical, the function $\Delta(1,2)$ as defined in Equation \eqref{eq:Delta_ij} is zero, \textit{i.e.}
		\[\det\begin{pmatrix}
		1			& \tilde{x}			& \tilde{y} 	& \tilde{z} \\
		1 & x		& y & z \\
		0 & 1 	& g' & u_x+g'u_y\\
		0 & 0 & g''		& u^{(2)}
		\end{pmatrix}=(x-\tilde{x})[g'u^{(2)}-g''(u_x+g'u_y)]-(y-\tilde{y})u^{(2)}+(z-\tilde{z})g''=0\]
		where $u^{(2)}= u_{xx}+2g'u_{xy}+g''u_y+(g')^2u_{yy}$.
		Substituting the equation
		\begin{equation}\label{eq:frobTPS}
			z-\tilde{z}=u_x(x-\tilde{x})+u_y(y-\tilde{y})
		\end{equation}
		 for the condition $\Phi(P) \in T_P S$ gives
		\[
		(x-\tilde{x})(g'u^{(2)}-g''g'u_y)-(y-\tilde{y})(u^{(2)}-g''u_y) =\left[(x-\tilde{x})g'-(y-\tilde{y})\right]\left(u^{(2)}-g''u_y\right)=0.
		\]
		
		Suppose that $g'=(y-\tilde{y})/(x-\tilde{x})$. Let us prove that then $\nu_1 > 1$ by noticing that the functions $\Delta(1,j)$ defined in Equation \eqref{eq:Delta_ij} are all zero for $j \geq 2$. Expanding the determinant $\Delta(1,j)$ and using the relation \eqref{eq:frobTPS} as we did with $\Delta(1,2)$, we get 
		\[\Delta(1,j)= \left[(x-\tilde{x})g'-(y-\tilde{y})\right]\left(u^{(j)}-g^{(j)}u_y\right)=0, \]
		as the first factor is zero by hypothesis. Since this contradicts Proposition \ref{propSV}, we conclude that $\nu_1>1$.
		
		\medskip
		
		Suppose now that $u^{(2)}-g''u_y=0$, that is 
		\begin{equation}\label{secondcase}
		u_{yy}(g')^2+2g'u_{xy}+u_{xx}=0.
		\end{equation}
		Deriving the equation $u-u(\tilde{x},g(\tilde{x}))=u_x(x-\tilde{x})+u_y(y-\tilde{y})$ we get
		\begin{equation}\label{secondderivative}
		u_{xx}(x-\tilde{x})+u_{xy}[g'(x-\tilde{x})+(y-\tilde{y})]+u_{yy}g'(y-\tilde{y})=0.
		\end{equation}
		Replacing $u_{xx}=-u_{yy}(g')^2-2g'u_{xy}$ in Equation \eqref{secondderivative} gives
		\[
		u_{yy}(\tilde{x}-x)(g')^2+g'[u_{xy}(\tilde{x}-x)+u_{yy}(y-\tilde{y})]+u_{xy}(y-\tilde{y})=0.
		\]
		The solutions of this quadratic equation in the variable $g'$ are
		$g'=(y-\tilde{y})/(x-\tilde{x})$ and $g'=-u_{xy}/u_{yy}$. 
		We have already dealt with the first case, so let us assume that $g'=-u_{xy}/u_{yy}$. 
		Since $u_{xx}+2u_{xy}g'+u_{yy}(g')^2=0$, we have
	\[u''=u_{xx}+2u_{xy}g'+u_{yy}(g')^2+g''u_y=g''u_y\]
	and $u^{(3)}=g^{(3)} u_y + g''(u_{xy}+g' u_{yy})=g^{(3)} u_y$. One can easily deduce that
	\begin{equation}\label{eq:derive}
		u^{(j)}=g^{(j)}u_y \text{ for every }j \geq 2.
	\end{equation}
	Now we conclude using Lemma \ref{lem:delta}.
	\end{proof}

We remark that if $S$ is an irreducible Frobenius non--classical surface, then $\CPhi^S$ is the whole surface $S$. Thus, the previous proposition entails that non--degenerate Frobenius non--classical curves with $\nu_1=1$ cannot lie on a Frobenius non--classical surface.
\subsubsection{Frobenius non--classical curves with $\nu_1>1$.}\label{ss:FNC}

Frobenius non--classical curves $C\subset S$ with $\nu_1>1$ are components of $\CPhi^S$, as the $q$--deleted order being $\epsilon_1=1$ implies that a general point $P$ on such curves satisfies $\Phi(P)\in T_P^{(1)} C \subset T_P S$ (see Equation \eqref{eq:deleted_order}). However, we will show that curves on $S$ with $\nu_1>1$ and having degree lower than or equal to the size of the finite field $\Fq$ are in fact degenerate components of $\CPhi^S$ (Proposition \ref{cor:nu1>1}). It means that these curves may have more points than prescribed by our bound, as proven by Borges and Homma \cite[Theorem 1.3]{BH17}.

\begin{lemma}\label{lem:pos_order_seq}
	Assume that $C$ is Frobenius non--classical with $\nu_1>1$ and $\delta \leq q$. Then $\nu_1=p^e$ and $\nu_2=p^{e'}$ or $\nu_2=2p^{e}$, for some integers $e < e'$.
\end{lemma}

\begin{proof}
	By the work of Homma \cite{H91} the possible order sequences $\{\epsilon_0,\epsilon_1,\epsilon_2,\epsilon_3\}$ of $C$ are
	\begin{enumerate}
		\item $\{0,1,2,3\}$ for $p>3$,
		\item $\{0,1,2,p^e\}$ for $p>2$,
		\item\label{it:NR2} $\{0,1,p^e,2p^e\}$,
		\item\label{it:NN} $\{0,1,p^e,p^{e'}\}$ for $p>2$,
		\item $\{0,1,p^e,p^e+1\}$,
	\end{enumerate}
	with integers $e < e'$. With the assumption $\nu_1 > 1$, we know that the $q$--deleted order is $\epsilon_1=1$, hence $\nu_1=\epsilon_2$ and $\nu_2=\epsilon_3$. 
	
	  Let us prove that the assumption on the degree of $C$ entails that only the sequences \eqref{it:NR2} and \eqref{it:NN} can occur. 
	By \cite[Corollary 5]{GH}, if $\delta \leq q$ then the $q$--deleted order, here $\epsilon_1=1$, is a $p$--maximal order, \textit{i.e.} if $\binom{\epsilon_j}{\epsilon_1} \neq 0 \mod p$, then $j=1$ (see \cite[p.31]{GH} for details). In the first two cases, $\binom{\epsilon_2}{\epsilon_1}=2 \neq 0 \mod p$. In the last case, we have $\binom{\epsilon_3}{\epsilon_1} = 1 \mod p$. Therefore, $\epsilon_1=1$ is not a $p$--maximal order in these three cases, which concludes the proof.	
\end{proof}

\begin{proposition}\label{cor:nu1>1}
	Assume that $C$ is Frobenius non--classical with $\nu_1>1$ and of degree $\delta \leq q$. Then, $C$ is a degenerate curve.
\end{proposition}
	\begin{proof}
	Let $p$ denote the characteristic of the field. Let $j \geq 2$ and consider the function $\Delta(1,j)$ as defined in Equation \eqref{eq:Delta_ij}. Following the proof of Proposition \ref{prop:FNCA=Nice}, we get

	\[\Delta(1,j)= \left[(x-\tilde{x})g'-(y-\tilde{y})\right]\left(u^{(j)}-g^{(j)}u_y\right).\]
	Since $\nu_1>1$, this function must be zero.
	If the second factor is zero, then we can conclude by Lemma \ref{lem:delta} that $C$ is degenerate. Thus, let us assume that $g'=(\tilde y - y)/(\tilde x - x)$.  As $y=g(x)$ and $\tilde{y}=g(\tilde{x})$, taking the derivative on both side implies that $g''=0$. Then, for any $2\leq j< p$ we have $g^{(j)}=0$. Furthermore, since $u_{xx}+2u_{xy}g'+u_{yy}(g')^2=0$, we also have
	$u^{(2)}=u_{xx}+2u_{xy}g'+u_{yy}(g')^2+g''u_y=0$. Then, for any $2\leq j< p$ we get $u^{(j)}=0$. Regarding the derivative of order $p$ of $u$, we have
	\begin{equation}\label{eq:up}
	u^{(p)}=g^{(p)}u_y+g^{(p-1)}(u_{xy}+g'u_{yy})=g^{(p)}u_y.
	\end{equation}
	Our goal is to show that under our hypothesis, for every $j \geq p$, we have	$u^{(j)}=g^{(j)}u_y$. 
		
	Consider the function $\Delta(i,j)$ for $i,j\geq 2$, that is
\[\Delta(i,j)=(x -\tilde x)\left(u^{(j)}g^{(i)}-u^{(i)}g^{(j)}\right).\]
We remark that if $(x -\tilde x) =0$ for any point of the curve, then $\Delta(i,j)=0$ which contradicts Proposition \ref{propSV}. Thus, we will suppose $(x -\tilde x) \neq 0$ in what follows.
	
	By Lemma \ref{lem:pos_order_seq}, we have $\nu_1=p^e$ and $\nu_2=p^{e'}$ or $\nu_2=2p^{e}$, for some integers $e < e'$.
	
	Assume that $\nu_1=p$, that is $e=1$. Since $\nu_2 > p+1$ we have
	\begin{equation}\label{eq: Deltapp+1}
		\Delta(p,p+1)=(x-\tilde{x}) \left(u^{(p+1)}g^{(p)}-u^{(p)}g^{(p+1)}\right)=0.
	\end{equation}
	Notice that  $g^{(p)}\neq 0$. Otherwise we would have $u^{(p)}=0$ and thus $\Delta(p,j)=0$ for any $j>p$, again contradicting Proposition \ref{propSV}. Thus Equation \eqref{eq: Deltapp+1} gives
	\[u^{(p+1)}=g^{(p+1)}\frac{u^{(p)}}{g^{(p)}}=g^{(p+1)}u_y,\] 
	where the last equality holds using Equation \eqref{eq:up}. Deriving Equation \eqref{eq:up} we have
	\[u^{(p+1)}=g^{(p+1)}u_y+g^{(p)}(u_{xy}+g'u_{yy}).\]
	Comparing the above expressions of $u^{(p+1)}$, we deduce that $g'=-u_{xy}/u_{yy}$. As in the end of Proposition \ref{prop:FNCA=Nice} this entails that $u^{(j)}=g^{(j)}u_y$ for any $j\geq 2$ and we can conclude that $C$ is degenerate using Lemma \ref{lem:delta}.
		
	Now, let us assume $\nu_1=p^e$ with $e>1$. Then for any $j>p$ we have
	\[\Delta(p,j)=(x-\tilde{x})\left(u^{(j)}g^{(p)}-u^{(p)}g^{(j)}\right)=0.\]
Thus for any $j>p$ we have
	\[u^{(j)}=g^{(j)}\frac{u^{(p)}}{g^{(p)}}=g^{(j)}u_y.\]
	Again by Lemma \ref{lem:delta} we conclude that $C$ is degenerate.
		\end{proof}

We remark that the previous proposition shows that curves with $\nu_1>1$ in $\PP^3$ of degree smaller than or equal to the size of the field are necessarily plane curves. Let us emphazise that this also holds when $S$ is Frobenius non--classical. Also, it is worth noting that one can find Frobenius non--classical curves with $\nu_1>1$ that are non--degenerate if their degree is larger than $q$, as shown in the following example.

	\begin{example}\label{ex:non-deg}
	 Over $\F_2$, on the singular quartic surface defined by
	 \begin{align*}
	f=& X_0^2X_1^2 + X_0X_1^3 + X_0^3X_2 + X_0^2X_1X_2 + X_1^3X_2 + X_0^2X_2^2 + X_1^2X_2^2 + X_0X_2^3 \\
	&+ X_1X_2^3 + X_2^4 +	X_0^3X_3 + X_1^3X_3 + X_0X_1X_2X_3 + X_2^3X_3 +	X_0^2X_3^2 \\
	&+ X_1^2X_3^2 + X_0X_2X_3^2 + X_1X_2X_3^2 + 	X_2^2X_3^2 + X_0X_3^3 + X_1X_3^3,
 	\end{align*}
 the curve $\CPhi^S$ splits into some rational lines, one singular plane curve of degree $3$ with $4$ rational points, and two non--degenerate curves: a non--singular quintic and a singular sextic, both having $5$ $\F_2$--points. Using {\scshape Magma}, we checked that all the components of $\CPhi^S$ have $\nu_1>1$.
	\end{example}

So far we have proved that Frobenius non--classical curves with $\nu_1>1$ and $\delta\leq q$ are degenerate components of $\CPhi^S$. When the surface is smooth we can completely characterize the degenerate components of $\CPhi^S$ as follows.

	\begin{proposition}\label{prop:non-sing-S}
		If the surface $S$ is smooth the degenerate non--linear components of $\CPhi^S$ have $\nu_1 >1$.
	\end{proposition}
	\begin{proof}
		Take $C \subset S$ a degenerate non--linear component of $\CPhi^S$. Denote by $H$ the plane containing $C$. If $H\neq T_P S$ for a general point, then $\Phi(P)\in H\cap T_P S = T_P^{(1)}  C$, thus $\nu_1>1$. Otherwise, suppose that $H=T_P S$ for all points $P\in C$. Then, the Gauss map 
		\[\gamma :\left\{\begin{array}{ccc}
			S& \rightarrow& (\PP^3)^*\\
			P & \mapsto & T_P S
		\end{array}\right.\]
		contracts the curve $C$. This contradicts the finiteness of the Gauss map of a smooth surface \cite[Corollary I.2.8]{Z93}.
		
	\end{proof}
	It is worth noting that this result together with Proposition \ref{prop:FNCA=Nice} shows that Frobenius non--classical curves with $\nu_1=1$ cannot lie on \emph{smooth} Frobenius non--classical surfaces.

	\subsubsection{Frobenius classical curves}\label{ss:FC}
	 To handle Frobenius classical components of $\CPhi^S$ we ran numerical implementations in {\scshape Magma}. They indicate that, when the curve $\CPhi^S$ (of degree $d(d+q-1)$) is reducible, its Frobenius classical non--degenerate components have degree strictly bigger than the size of the finite field. In particular, we checked all the quadrics over $\F_2$ and $\F_3$. We conjecture this is true in general.
	
	
		\begin{conj}\label{conjecture}
		Non--degenerate Frobenius classical components of the curve $\CPhi^S$ have degree larger than $q$.
	\end{conj}

	The condition on the degree is tight, as illustrated by the following examples we studied using {\scshape Magma}.
	
	\begin{example}\label{ex1}
		Along our tests on cubic surfaces, the non--singular surface we encountered on which the curve $\CPhi^S$ has the highest number of components is defined over $\F_5$ by 
		\begin{align*}
			f=& 2X_0X_1^2 + 2X_1^3 + 2X_0^2X_2 + 2X_0X_1X_2 + X_1^2X_2 + 2X_0X_2^2 + 3X_1X_2^2 + 3X_2^3\\
			& -X_0^2X_3+ X_0X_1X_3 + X_1^2X_3 + 2X_1X_2X_3 + 2X_2^2X_3 +	3X_0X_3^2 -X_1X_3^2 + X_2X_3^2.
		\end{align*}
		The curve $\CPhi^S$ has degree 21 and is formed of 15 $\F_5$--lines and one non--degenerate sextic $C$. The curve $C$ is Frobenius classical and it has genus $4$ and it contains $18$ $\F_5$--points.
		
		The number of rational points on this sextic curve still satisfies the bound \eqref{notreborne}. Indeed, the curve $C$ is defined by $f=f^{(1)}=f^{(2)}=0$ where
		\begin{align*}
		f^{(1)}&=X_0^2 + 2X_0X_1 + 3X_1^2 + 2X_0X_2 + X_1X_2 + 3X_0X_3 +
			4X_2X_3 + X_3^2,\\
		f^{(2)}&=X_0X_2X_3 + X_2^2X_3 + 2X_0X_3^2 + X_1X_3^2 +
		X_2X_3^2 + 3X_3^3.
		\end{align*}
		The quadric surface $S_1$ defined by $f^{(1)}=0$ is Frobenius classical by Lemma \ref{FCS}. We checked that the curve $C$ is not a component of the curve $\CPhi^{S_1}$. Thus we can apply our bound for $C\subset S_1$, which gives $\#C(\Fq) \leq \frac{1}{2}\cdot 6(2+5-1)=18$.
		\end{example}
		
	\begin{example}\label{ex2} Over $\F_3$, on the \emph{singular} cubic surface defined by $f=0$, where
		\begin{align*}		
		f=&-X_0^3 -X_1^3 + X_0^2X_2 + X_0X_1X_2 + X_1^2X_2 + X_0X_2^2 -X_1X_2^2 -X_0X_1X_3 \\
		&+ X_1^2X_3 + X_0X_2X_3 -X_2^2X_3 + X_0X_3^2 + X_1X_3^2 -X_2X_3^2 + X_3^3,
		\end{align*}
		the degree--15 curve $\CPhi^S$ splits into one $\F_3$--line (counted three times) and three non--degenerate quartics, which are Frobenius classical. These curves are all isomorphic, of genus $1$, and with $7$ $\F_3$--points. As in the previous example, they lie on a Frobenius classical quadric surface $S_1$, on which they do not belong to the curve $\CPhi^{S_1}$. Hence, they fall under our bound for $(d,\delta)=(2,4)$.
		
	\end{example}

	\subsection{Impact on the main theorem} 
	So far we have characterised the irreducible components of $\CPhi^S$ in light of Theorem \ref{th:P3}. On the one hand, since non--degenerate Frobenius non--classical curves on $S$ with $\nu_1=1$ are never components of $\CPhi^S$, the bound \eqref{notreborne} always applies to these curves. On the other hand, Frobenius non--classical curves with $\nu_1>1$ and degree $\delta\leq q$ are plane curves, so they are beyond the scope of our paper. Concerning Frobenius classical curves, if Conjecture \ref{conjecture} holds, then they are not components of $\CPhi^S$ whenever their degree is smaller than or equal to $q$. We sum up these observations in the following corollary.
	\begin{corollary}\label{cor:P3}
		Let $S$ be an irreducible normal Frobenius classical surface defined over $\Fq$ of degree $d>1$ in $\PP^3$. Let $C$ be a non--degenerate irreducible curve defined over $\Fq$ of degree $2< \delta \leq q$ lying on $S$. Suppose $C$ is Frobenius non--classical. Then
		\[\#C(\F_q)\leq \frac{\delta(d+q-1)}{2}. \]
		Moreover, under Conjecture \ref{conjecture}, the bound holds for Frobenius classical curves too.
	\end{corollary}

%

\section{Comparison and conclusions}\label{Sec:conclusion}

	\subsection{Comparison}
In this subsection, we compare the bound from Theorem \ref{th:P3} to other existing bounds for irreducible projective space curves.
\subsubsection{On smooth quadrics}

Couvreur and Duursma computed the minimum distance of evaluation codes on smooth quadrics \cite{CD13}, which gives an upperbound for the number of rational points on curves lying on these surfaces. There are two classes of smooth quadrics up to $\Fq$--isomorphism in $\PP^3$: the hyperbolic quadric $\mathcal{H}$, isomorphic to the image of $\PP^1 \times \PP^1$ under the Segre embedding, and the elliptic twisted quadric $\mathcal{E}$. Couvreur and Duursma's upper bound for curves obtained by intersecting $\mathcal{H}$ with a surface of degree $s$ is reached only by union of lines \cite[Corollary 2.5]{CD13}. On the twisted quadric $\mathcal{E}$, their bound equals $s(q+1)$ \cite[Corollary 3.13]{CD13} and matches the one given by Theorem \ref{th:P3} ($d=2$ and $\delta=2s$). They proved it may be reached by irreducible curves \cite[Example 3.14]{CD13}.

\subsubsection{Numerical comparison to other bounds} 

Figure \ref{fig:comp} illustrates how Theorem \ref{th:P3} improves preexisting bounds on the number of $\Fq$--points on space curves \cite{SV86,H12} for some range of parameters.

Let us first notice that if the degree $d$ of the surface $S$ is larger than or equal to $q$, then Homma's bound for non--degenerate space curves \cite{H12}, that is $\# C(\Fq) \leq q(\delta-1)+1$, is sharper than the bound provided by Theorem \ref{th:P3}.

St\"ohr and Voloch \cite{SV86} proved that a non--singular space curve $C$ defined over $\Fq$ of degree $\delta$, genus $g$ and with Frobenius orders $\nu_1$ and $\nu_2$ satisfies
\[\#C(\Fq) \leq \frac{1}{3}\left(2(\nu_1+\nu_2)(g-1)+ \delta(q+3)\right). \]
To compare Theorem \ref{th:P3} to St\"ohr and Voloch's formula, we use Harris's result \cite{H80} to bound the genus from above, and we assume that the curve is Frobenius classical, \textit{i.e.} $\nu_1=1$ and $\nu_2=2$. 


\tikzset{
	declare function={Floor(\x)=round(\x-0.5);},
}

\begin{figure}[h]
	\centering
	\begin{subfigure}[t]{.45\linewidth}
		\centering
		\begin{tikzpicture}[declare function={ 
				eps(\x,\d)=round(-\x/\d-Floor(-\x/\d));
				PI(\x,\d) = \x/2*(\x/\d +\d-4)+1-eps(\x,\d)/2*(\d-eps(\x,\d)-1+eps(\x,\d)/\d) ;
				g(\x,\d) = 	(\x>\d*(\d-1)) * PI(\x,\d)+ (\x<=\d*(\d-1)) * PI(\x,Floor((\x-1)/\d)+1);},]
			\def\q{9};
			\def\d{5};
			\begin{axis}[xmin=0,xmax=\q,ymin=0,ymax=80,clip=true,axis lines*=middle, grid = both,
				grid style={dashed}, xtick distance=2, extra x ticks={0,1,...,\q},extra x tick labels={},ytick distance=20,x=0.5cm,legend pos=south east, xlabel={Degree $\delta$ of the curve $C$}, ylabel={\small Upper bound},  domain={0:\q},samples=\q+1]
				\addplot[black,thick] {x*(\d+ \q-1)/2};	
				\addplot[very thick,dotted] {\q*(x-1)+1}; 
				\addplot[densely dashed,very thick] {4*(g(x,\d)-1)+(\q+3)*x/3}; 
			\end{axis}
		\end{tikzpicture} 
		\subcaption{$q=9$ and $d=5$}
	\end{subfigure}
	\qquad
		\begin{subfigure}[t]{.45\linewidth}
		\centering
		\begin{tikzpicture}[declare function={ 
				eps(\x,\d)=round(-\x/\d-Floor(-\x/\d));
				PI(\x,\d) = \x/2*(\x/\d +\d-4)+1-eps(\x,\d)/2*(\d-eps(\x,\d)-1+eps(\x,\d)/\d) ;
				g(\x,\d) = 	(\x>\d*(\d-1)) * PI(\x,\d)+ (\x<=\d*(\d-1)) * PI(\x,Floor((\x-1)/\d)+1);},]
			\def\q{13};
			\def\d{4};
			\begin{axis}[xmin=0,xmax=\q,ymin=0,ymax=120,clip=true,axis lines*=middle, grid = both,
				grid style={dashed}, xtick distance=2, extra x ticks={0,1,...,\q},extra x tick labels={},ytick distance=30,x=0.4cm,legend pos=south east, xlabel={Degree $\delta$ of the curve $C$},  domain={0:\q},samples=\q+1]
				\addplot[black,thick] {x*(\d+ \q-1)/2};	
				\addplot[very thick,dotted] {\q*(x-1)+1};
				\addplot[densely dashed,very thick] {4*(g(x,\d)-1)+(\q+3)*x/3}; 
			\end{axis}
		\end{tikzpicture} 
		\subcaption{$q=13$ and $d=4$}
	\end{subfigure}
	\\
	\begin{subfigure}[b]{\textwidth}
		\centering
		\begin{tikzpicture} 
			\begin{axis}[%
				hide axis,
				xmin=10,
				xmax=50,
				ymin=0,
				ymax=0.4,
				legend style={draw=white!15!black,legend cell align=left},legend columns=-1,
				/tikz/every even column/.append style={column sep=0.2cm}]
				\addlegendimage{black} \addlegendentry{Theorem \ref{th:P3}};
				
				\addlegendimage{very thick,dotted} \addlegendentry{\cite{H12}};
				

				\addlegendimage{densely dashed} \addlegendentry{\cite{SV86}};

			\end{axis}
		\end{tikzpicture}.

	\end{subfigure}
	\vspace*{-5cm}
	\caption{Comparison of bounds on the number of $\F_{q}$--points on a non--degenerate curve $C$ of degree $\delta$ lying on a surface $S \subset \PP^3$ of degree $d$.}\label{fig:comp}
\end{figure}
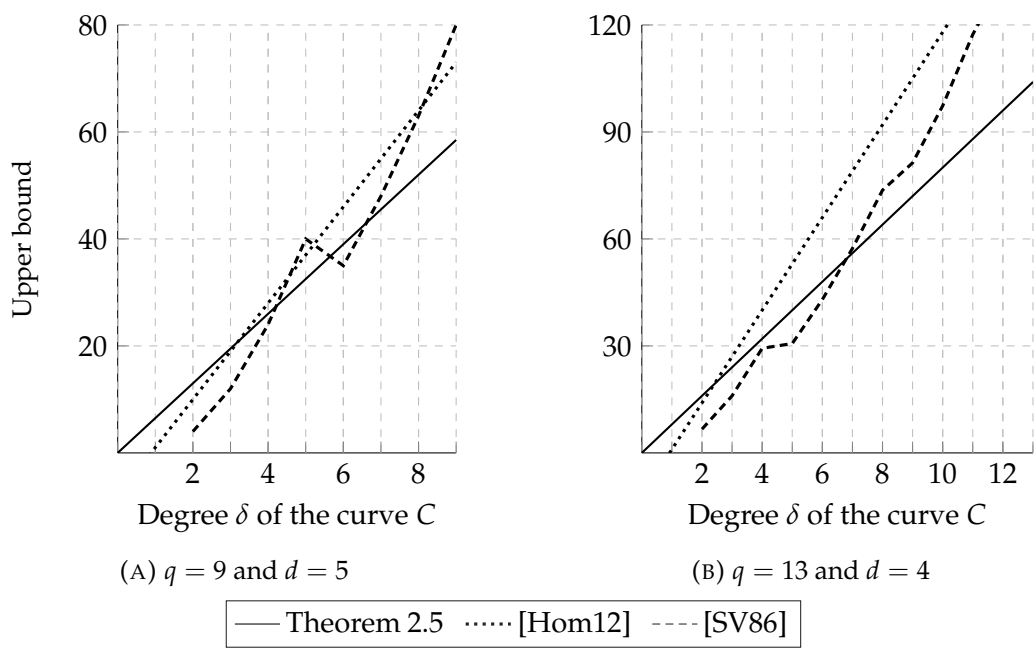

\begin{example}
	For a curve of degree $11$ on a degree $5$ surface defined over $\F_9$, the bound of Theorem \ref{th:P3} equals 72. Such a curve has genus at most $17$. Curves of genus $16$ and $17$ referenced in the database of curves with many points \cite{ManyPoints} over $\F_9$ have $74$ $\F_9$--points. In particular, this means that record curves cannot be found on a Frobenius classical surface $S$ for this choice of degrees, unless being components of $\CPhi^S$. This illustrates how embedding entails arithmetic and geometric constraints on a variety.
\end{example}

\subsection{Future works}

Theorem \ref{th:P3} essentially relies on the intersection theory in $\PP^3$. It is a natural question to ask whether our approach can be generalised to bound the number of rational points on curves lying on surfaces embedded in $\PP^n$, for $n>3$.

Let us assume that $S$ is an irreducible algebraic surface defined over a finite field $\Fq$ and embedded in a projective space $\PP^n$, for $n \geq 3$. Let $C$ be a curve lying on $S$. In $S \times \PP^n$, let us consider the two varieties defined as follows:
\begin{itemize}
	\item $\Gamma_C= \{(P,\Phi(P)) \in C^2 \mid P \in C\}$ the graph of $\Phi$ restricted to the curve $C$,
	\item $\mathcal{T}_S=\{(P,Q) \in S \times \PP^n \mid P \in S, \: Q \in T_P S\}.$
\end{itemize}

The intersection of the varieties $\Gamma_C$ and $\mathcal{T}_S$ is in one--to--one correspondance with the points $P$ on $C$ whose image under the Frobenius map lies on their own tangent plane $T_P S$. Of course, the rational points on $C$ can be embedded into this intersection via the diagonal morphism on $C$.  In fact, the set $\CPhi^S$ we study in this paper is the image of $\Gamma_C \cap \mathcal{T}_S \in S\times \PP^3$ under the first projection.

The variety $\Gamma_C$ is isomorphic to $C$ and it is thus a curve of $S \times \PP^n$. The variety $\mathcal{T}_S$ has a dimension $4$ if $S$ is non--singular. Consequently, these two varieties have complementary dimensions in $S\times \PP^n$ (of dimension $n+2$) if and only if $n=3$. In this paper, we have bounded the number of rational points on $C$ by a fraction of the intersection product $[\Gamma_C]\cdot [\mathcal{T}_S]$. When $n \geq 4$, this intersection product is zero, whereas the intersection $\Gamma_C \cap \mathcal{T}_S$ is nonempty.  A common way to fix this dimension incompatibility is by blowing up one of the varieties, for instance $\mathcal{T}_S$ or $S \times S$. We are likely to run into considerable difficulties in this context but we hope to return to it in a later paper.

\bigskip

\noindent
\textbf{Acknowledgements.} 
The two authors are grateful to J. F. Voloch for his thoughtful hints. They also thank the anonymous referee whose comments and questions drove them to improve the paper.\\
This work was funded in part by the grant ANR-21-CE39-0009-BARRACUDA from the French National Research Agency. The first author has also received funding from the French “Agence de l'innovation de défense”, from the French ‘‘Fondation Mathématique Jacques Hadamard", and  from  the  European Union’s Horizon 2020 research and innovation programme under the Marie Skłodowska--Curie grant agreement No 899987. The second author is supported by the French government ‘‘Investissements d’Avenir" program ANR-11-LABX-0020-01.

	\bibliography{biblioCSPn}
	\bibliographystyle{alpha}

\end{document}